\documentclass[11pt,reqno]{amsart}
\usepackage{amssymb,verbatim}
\usepackage[hypertex]{hyperref}
\usepackage[notref,notcite]{showkeys}
\makeatletter
\def\part{\@startsection{part}{0}%
\z@{\linespacing\@plus\linespacing}{.5\linespacing}%
{\normalfont\scshape\centering}}                % see: section
\def\section{\@startsection{section}{1}%
\z@{.7\linespacing\@plus\linespacing}{.5\linespacing}%
{\normalfont\bfseries\centering}}    % <=== (цембфемшоп рпддетцйчбфш)  FONT : BF-SERIES IN THE SECTION HEADING
\def\l@part{\@tocline{1}{6pt plus 1pt}{0pc}{}{\hfil\scshape}} % vskip 6pt plus 1pt
\def\l@section{\@tocline{1}{4pt}{0pc}{}{}} % vskip 4 pt
\def\l@subsection{\@tocline{2}{0pt}{0.4pc}{5pc}{}}
\def\l@subsubsection{\@tocline{3}{0pt}{0.4pc}{7pc}{}}
\makeatother
\makeatletter
\def\@setaddresses{\par
  \nobreak \begingroup
\footnotesize
  \def\author##1{\nobreak\addvspace\bigskipamount}%
  \def\\{\unskip, \ignorespaces}%
  \interlinepenalty\@M
  \def\address##1##2{\begingroup
    \par\addvspace\bigskipamount\indent
    \@ifnotempty{##1}{(\ignorespaces##1\unskip) }%
    {\scshape\ignorespaces##2}\par\endgroup}%
  \def\curraddr##1##2{\begingroup
    \@ifnotempty{##2}{\nobreak\indent{\itshape Current address}%
      \@ifnotempty{##1}{, \ignorespaces##1\unskip}\/:\space
      ##2\par}\endgroup}%
  \def\email##1##2{\begingroup
    \@ifnotempty{##2}{\nobreak\indent{\itshape \EMAILADDRESS}%
      \@ifnotempty{##1}{, \ignorespaces##1\unskip}\/:\space
      \rmfamily##2\par}\endgroup}%
  \def\urladdr##1##2{\begingroup
    \@ifnotempty{##2}{\nobreak\indent{\itshape URL}%
      \@ifnotempty{##1}{, \ignorespaces##1\unskip}\/:\space
      \ttfamily##2\par}\endgroup}%
  \addresses
  \endgroup
}
\makeatother
\def\EMAILADDRESS{E-mail address}
  \let\phi\varphi
\let\opn\operatorname \let\opl\operatornamewithlimits
\def\vac{{1\!\!\mathbb I}}

\theoremstyle{remark}\newtheorem*{REM*}{Remark}
\theoremstyle{plain}
\newtheorem*{LEM}{Lemma}%[section]
\newtheorem{THM}{Theorem}%[section]

\address{A.~M.~Vershik, St.~Petersburg Department of Steklov Institute of
Mathematics, 27 Fontanka, 191023 St.~Petersburg, Russia}
\email{vershik@pdmi.ras.ru}
\address{M.~I.~Graev, Institute for System Studies, 36-1
Nakhimovsky pr., 117218 Moscow, Russia}
\email{graev\_36@mtu-net.ru}
\date{August 9, 2008}
\title[\today]{}
\keywords{%лМАЮЕЧЩЕ УМПЧБ.
Current group, integral model, Fock representation, special representation,
infinite-dimensional Lebesgue measure}
\begin{document}
%__LUDE_

{\tiny\par\noindent\today \hfill newimod.tex 09.08.2008}

\noindent
UDC 517.5

\begin{center}{}\bf
Integral models of unitary representations of current groups
with values in semidirect products
\\ \rm
A.~M.~Vershik \footnote{Supported by the grants NSh-2460.2008.1
and RFBR 08-01-00379.},
M.~I.~Graev \footnote{Supported by the RFBR grant 07-01-00101a.}
\end{center}

\begin{abstract}
{We describe a general construction of irreducible unitary representations
of the group of currents with values in the semidirect product of a locally compact
subgroup $P_0$ and a one-parameter group ${\mathbb R {\,}}^*_+=\{r:r>0\}$
of automorphisms of $P_0$. This construction is determined by a
a faithful unitary representation of $P_0$ (canonical representation)
whose images under the action of
the group of automorphisms tend to the identity representation as
$r\to 0$. We apply this construction to the groups of currents of the maximal
parabolic subgroups of the groups of motions of the $n$-dimensional
real and complex Lobachevsky spaces. The obtained representations
of the groups of parabolic currents can be uniquely extended to the groups
of currents with values in the semisimple groups
$O(n,1)$ and $U(n,1)$. This gives a new description of the representations of
the groups of currents of these groups constructed in the 70s
and realized in the Fock space. The key role in our construction is played
by the so-called special representation of the parabolic subgroup $P$
and the remarkable $\sigma$-finite measure (Lebesgue measure)
$\mathcal L$ in the space of distributions.}

\end{abstract}
\vskip-1cm
\maketitle
\vskip-0.5cm\vskip-0.5cm%\vskip-0.5cm%\vskip-0.5cm

\begin{flushright}{}
\sl
To dear Israel Moiseevich Gelfand on the occasion of his 95th birthday
\end{flushright}

\bigskip

In 1972, on the initiative of I.~M.~Gelfand, a series of papers by three authors
(I.M., M.I., A.M) on unitary representations of functional groups, or
current groups, was started. The problem that was posed by I.M. to the first author of
this paper in spring 1972 and that initiated this series of papers
was to find out whether there exists a ``multiplicative integral
of representations'' (see below) for the group $SL(2,R)$.
Soon it became clear that the answer to this question is positive,
and the first paper of the series \cite{VGG-1} contained a description
of several models of this representation. The main idea was to study a
neighborhood of the identity representation of the group
$SL(2,R)$ itself and the so-called special (infinitesimal) representation
of this group, whose first cohomology is nontrivial. A multiplicative integral
of representations is, in the authors' words, a ``tensor product
of infinitely many infinitesimal representations.'' In the subsequent joint
papers, various generalizations of this construction (to simple Lie
groups of rank~$1$, for groups of diffeomorphisms, etc.) were found,
and techniques for working with such representations were developed.
Starting from 2004, the authors of this paper initiated a systematic study
of integral models of representations of current groups based on a new
interpretation of continual tensor product that is different from the Fock one
and essentially exploits a remarkable $\sigma$-finite measure in the space
of distributions. A systematic presentation of the whole area
will be given in the book {\it The Representation Theory of Current Groups}
by the three authors, which is now in preparation.

\section{Introduction}
We give a brief description of the theory of integral models of representations
of current groups introduced in the previous paper by the same authors
\cite{VG-7}. The integral models for the groups of currents of the parabolic
subgroups of semisimple groups of  rank~$1$ are used for constructing representations
of the groups of currents of these semisimple groups themselves.
We also establish a relation between the integral models and
the corresponding representations
of the groups of currents of semisimple groups of rank~$1$
and the Fock models considered earlier
(\cite{VGG-6,GG}).

By a group of currents we mean a group of measurable functions, on some space or
manifold equipped with a measure,
with values in some group of coefficients, endowed with the pointwise multiplication.
Only those irreducible representations of current groups are of interest
that are invariant under some group of automorphisms of the base space, for
example, the group of rotations of the circle (in the case of loops)
or the group of all measure-preserving transformations.
In \cite{VGG-1}, such a representation (if it does exist) was called
 a {\it multiplicative integral of representations}, since
in a sense it is a continuous tensor product of representations of the
group of coefficients. The problem of existence of such
a representation for a given group of coefficients and a given base space
is far from being trivial. The general scheme of studying representations
of current groups in a {\it Fock space} (i.e., in a Hilbert space
endowed with a Fock factorization respected by representations)
is presented in \cite{Araki}. In particular, in this case
the base of the current group
can be an arbitrary measure space, and the group of coefficients $G$ must have
an irreducible unitary representation with nontrivial first cohomology.
As shown \cite{VK}, the latter condition
can be satisfied only if the identity representation of $G$ is not isolated
in the set of irreducible unitary representations with the Fell topology
(see \cite{XR}). Using the adopted terminology, one can say that $G$
must not satisfy Kazhdan's property~(T). An irreducible unitary representation
$\pi$ with nontrivial first cohomology $H^1(G;\pi)$ is called
{\it special}. Among the semisimple Lie groups,
only some groups of rank~$1$ have  a special representation,
namely, the groups of motions of
the real and complex Lobachevsky spaces:
$U(n,1)$ and $O(n,1)$. For the groups of currents with coefficients
in these groups, the
corresponding irreducible representations  were first
constructed in \cite{VGG-1,VGG-2}. They are realized in
a Fock space, which is the exponential of a one-particle space, the latter
being in turn a direct integral of special representations. The construction
of any representation
of a current group in a Fock space
substantially uses a nontrivial $1$-cocycle of the group of coefficients
with values in a special representation of this group, whose existence guarantees the irreducibility
of the representation of the current group. The Fock model is convenient,
and it is used in most papers on the representation theory
of infinite-dimensional groups. In this model, the operators of the maximal
compact subgroup of currents are diagonalized.
Note also that the structure of a Fock Hilbert space can be equivalently
described with the help of the infinite-dimensional Gaussian measure
(``white noise''), and that the model suggested below uses another remarkable
measure in an infinite-dimensional space.

The paper \cite{VGG-6} initiated the study of another model of representations
of current groups, in which the operators of unipotent subgroups are diagonalized.
This approach has gradually led to a quite new understanding of the nature
of multiplicative integrals of representations and continuous tensor products.
It is based on two principal observations.

1) As in the case of the classical semisimple Lie groups, e.g.,
$SL(2,\Bbb R)$, one should start constructing a representation of a group of semisimple
currents with constructing a {\it representation of currents of
the maximal parabolic subgroup}, and then extend it to the whole group
of currents of the semisimple group.

2) In turn, the existence of a representation
of the group of currents of
the maximal parabolic subgroup of a semisimple group of  rank~$1$ is due to the fact
that there exists a remarkable $\sigma$-finite measure on the cone of discrete finite
positive measures on a manifold that is {\it invariant under the continuous
analog of the Cartan subgroup}, namely, the group of multiplicators with
finite integral of the logarithm.
This measure, which was introduced in
\cite{VGG-6}, studied in
\cite{TVY,VT}, and investigated in detail in the recent papers
\cite{Vershik,V1}, was called the {\it infinite-dimensional Lebesgue
measure}. It allows us to drastically simplify all constructions of
representations of current groups of semidirect products. Moreover,
it proves to be closely related to many combinatorial, analytical,
and probabilistic problems.

For the reason that will be clear from the construction, we call this type
of models of representations of current groups integral models.
In this paper, we present a coherent treatment of this alternative
construction of {\it integral models of representations of
current groups of semidirect products}, and show how they can be extended
to the groups of currents of the groups of motions of the Lobachevsky spaces.

In brief, our model of continuous tensor products of Hilbert spaces
and representations of groups of currents (loops) in these continuous
tensor products is
as follows:

a) we consider the space of countable linear positive
combinations of $\delta$-measures on the base
space with a finite sum of coefficients,
endowed with the infinite-dimensional Lebesgue
measure defined below \cite{V1,Vershik};

b) with each trajectory we associate a countable tensor product of Hilbert spaces
and representations in these spaces, taken along this trajectory;

c) we define the direct integral of the constructed countable tensor products with respect to the Lebesgue measure
on the space of trajectories.

Our construction yields representations of groups of currents (and,
in particular, groups of loops) in a Hilbert space endowed with an additional
structure, denoted by $\opn{INT}$; these representations are equivalent
to the representations constructed earlier in
\cite{VGG-1, VGG-2} and realized in the Fock space
$\opn{EXP}$. The (unique up to a scalar) isomorphism between the spaces
$\opn{EXP}$ and
$\opn{INT}$ is of interest in itself.\footnote{Note that this isomorphism
preserves the structure of a {\it Fock factorization of tensor products}
of Hilbert spaces. At the same time, there exist non-Fock factorizations
in Hilbert spaces (see \cite{TsV}), and hence there is a possibility to
construct representations of current groups with factorizations nonisomorphic
to the Fock one. As far as we know, this possibility has not yet been studied.}
One more advantage of this method of constructing representations of current groups
is that the irreducibility of the representation of the
parabolic group of currents, as well as the analysis
of the special representation, becomes almost obvious. The only problem is
to extend representations
to the whole semisimple group of currents. This turns out to be possible for
$O(n,1)$ and $U(n,1)$, and impossible for $Sp(n,1)$.

In Section~2 we define the class of semidirect products that
serve as the groups of coefficients of current groups
in our construction
of integral models of representations.
We introduce the notion of a canonical representation, which is used
in the definition of a special representation of a semidirect product.
In Section~3 we describe the main construction. Sections~4--5 contain
formulas for the spherical function of the integral model and for the
isomorphism between the integral model and  the Fock model. In Section~6 we
give a detailed description of the
integral model in the simplest case, that of the parabolic subgroup
of the group $ SL(2,{\mathbb R {\,}})$,
and construct an extension of this integral model to the whole group.
In Sections~7--8 we do the same for the general case
of the parabolic subgroups of the groups $U(n,1)$ and
$O(n,1)$.
For $U(n,1)$, we should distinguish the cases of orthogonal and unitary
special representations: in the unitary case, there arise projective
unitary representations of the group of currents of
$U(n,1)$, also considered in \cite{GG}. There is a relation between these
projective representations and unitary (nonprojective) representations of
the current group
$\widetilde G^X$, where $\widetilde G$ is the universal covering group of $U(n,1)$
constructed in \cite{B}.

Note also that there exists a well-developed theory of projective
highest weight representations  for the case of groups of smooth loops
with semisimple coefficients (the Kac--Moody theory).
The same applies also to other representations of smooth currents.
The proofs and the detailed version of the theory of integral models
is presented in another paper, presently under preparation.

\section{The original construction: special representations of semidirect
products of groups}

We begin with a general scheme which applies, in particular, to
parabolic subgroups of semisimple groups of rank~$1$. Consider the semidirect
product ${P}=S\rightthreetimes {P}_0$ of a locally compact group
$P_0$ and a one-parameter multiplicative group
$S\cong {\mathbb R {\,}}^*_+$ of automorphisms of $P_0$. Denote by $g^r$ the image
of an element $g\in P_0$ under an automorphism $r\in S$. In what follows,
unless otherwise stated, the term ``representation'' means an orthogonal
or unitary representation. The group of automorphisms also acts on
representations $T$ of the group $P_0$; namely, every representation $T$
of the subgroup ${P}_0$ in a Hilbert space
$H$ gives rise to a family of $T$-associated representations
$T_r$, $r\in {\mathbb R {\,}}^*_+$, obtained from $T$
by automorphisms from the group
${\mathbb R {\,}}^*$:  $T_r(g)\equiv T(g^r)$.

A cyclic representation $T$ of the subgroup $P_0$  is called {\it canonical}
with respect to the group of automorphisms $S$ if there exists a unit cyclic
vector $h\in H$ such that
$\|T_r(g)h-h\|<c(g)r$ for sufficiently small $r$ and any $g\in {P}_0$. A vector
$h$ satisfying this property is called almost invariant with respect to
the representation $T$. Obviously, in the Fell topology (see, e.g., \cite{XR})
on the space of representations, the representations
$\{T_r\}$ tend to the identity representation as $r\to 0$. In terms of
the matrix elements of $T$ (states) corresponding to cyclic vectors $h$,
our definition reads as follows:
$$\lim_{r \to 0}\frac{1-\phi_h(g^r)}{r}=-\frac{d\phi_h(g^r)}{dr}{\bigg|}_{r=0}\equiv  \psi(g)<\infty$$
for all $g\in G$, where $\phi_h(g)={\rm Re} \{\langle T(g)h,h\rangle\}$
is a normalized positive definite function on the group $G$
(a ``state'').\footnote{The function $\psi(\cdot)$ defined above is an unbounded
conditionally positive definite function and is equal to the squared norm
of the nontrivial cocycle in the special representation (see below).
As to the terminology, note that in \cite{VG-7} a canonical representation
was called summable, and in \cite{VGG-1} canonical states for the group $PSL(2,R)$
were defined as states that are invariant with respect to a compact subgroup
and satisfy a close condition. Here we somewhat extend the notion of
a canonical state using the same term.}

With every canonical representation $T$ of the subgroup
${P}_0$ we can associate a representation
$\widetilde T$ of the group ${P}={\mathbb R {\,}}^*_+\rightthreetimes {P}_0$
in the direct integral of the Hilbert spaces $H_r \equiv H$
with respect to the multiplicative Haar measure
$r^{-1}\,dr$ on
${\mathbb R {\,}}^*_+$,
$$
\mathcal H=\int_0^\infty H_r\,r^{-1}\,dr,
$$
i.e., in the space of sections $f(r)$ of the fiber bundle over
${\mathbb R{\,}}^*_+$ with fiber $H_r$ over a point $r\in{\mathbb R{\,}}^*_+$.
The action of the operators $\widetilde T(g_0)$, $g_0\in {P}_0$, on
$\mathcal H$ is fiberwise:
$(\widetilde T(g_0)f)(r)=T_r(g_0)f(r)$; the operators $\widetilde T(r_0)$,
$r_0\in {\mathbb R {\,}}^*_+$, permute the fibers:
$(\widetilde T(r_0)f)(r)=f(r_0r)$.
\begin{LEM}
The representation $\widetilde T$ of the group ${P}$ is a
special representation (see \cite{VG}).
\end{LEM}

Indeed, by construction, this representation has nontrivial $1$-cocycles
$b:\, {P}\to\mathcal H$ of the form
\begin{equation}{}\label{14-4} b(g,r)=(\widetilde
T(g)f_0)(r)-f_0(r),\quad f_0(r)=\phi(r)\,h_r,
\end{equation}
where $h_r=h$ is a vector in $H$ almost invariant with respect to $T$ and
$\phi(r)$ is an arbitrary smooth function on
$[0, \infty )$ such that
$\int_{\epsilon }^\infty \phi^2(r)r^{-1}dr < \infty $ for
$\epsilon >0$ and $\phi(0)\ne 0.$

Note that the $1$-cocycles corresponding to functions
$\phi$ and $\phi_1$ are equivalent if and only if
$\phi(0)=\phi_1(0)$.

Examples of groups of the form
${P}={\mathbb R{\,}}_+^*\rightthreetimes {P}_0$ for which
the subgroup $P_0$ has canonical representations
are the maximal parabolic subgroups of the groups
$O(n,1)$, $U(n,1)$, and $Sp(n,1)$. In the case of $O(n,1)$, there is only
one, up to passing to conjugate representations, canonical representation. In the cases of
$U(n,1)$ and $Sp(n,1)$, the group $P_0$ has countably many irreducible canonical
representations. Each of them gives rise to a special irreducible
representation of the group $P$.
Note that among all special representations of the group
$P \subset U(n,1)$, only two unitary representations and one real representation
can be extended to a special representation of the whole group
$U(n,1)$. The group $Sp(n,1)$  has no special representations at all.

\section{Construction of the integral model of representations ($\opn{INT}$)}

The {\it current group} ${P}^X$, where ${P}$ is a locally compact group
and $X$ is a space with a probability measure $m$, is
the group of bounded measurable mappings
$X\to {P}$ endowed with pointwise multiplication.
In \cite{VG-7}, for groups ${P}$ of the form
${P}={\mathbb R {\,}}^*_+\rightthreetimes {P}_0$, for each canonical representation
$T$ of the subgroup ${P}_0$ in a space $H$ and a vector
$h\in H$ almost invariant with respect to $T$, we defined a representation
of the current group ${P}^X$, which was called the integral model. By analogy
with Fock models $\opn{EXP}$ (see \cite{VGG-6}), we will denote this
representation of the group ${P}^X$ by $\opn{INT}\,T$, and
the space in which it is realized, by
$\opn{INT}\,H$, or, in more detail, $\opn{INT}\,(H,h)$.

The construction of the representation $\opn{INT}\,T$
of the current group ${P}^X$ for
${P}={\mathbb R {\,}}^*_+\rightthreetimes {P}_0$ is similar
to the construction of a special representation of the group $P$ from
a canonical representation $T$ of the subgroup
${P}_0$. First of all, in this construction, the space
${\mathbb R {\,}}^*_+$ is replaced by the cone of real positive
finite atomic measures on the space $(X,m)$:
$$
l^1_+(X)=\left\{\left. \xi=\sum_{k=1}^\infty  r_k \delta_{x_k}
\right| r_k>0,\, \sum_k r_k<\infty \right\};
$$
in what follows, elements of this space will be denoted by
${\xi}=\{r_k,x_k\}$. Instead of the multiplicative Haar measure
$r^{-1}dr$ on ${\mathbb R {\,}}^*_+$, we consider the
{\it infinite-dimensional $\sigma$-finite measure} $d{\mathcal L}(\xi)$ on
$l^1_+(X)$. This measure (see \cite{Vershik}), which is an
infinite-dimensional analog of the Lebesgue measure on the octant
${\mathbb R {\,}}^n_+$, is uniquely determined by its Laplace transform:
\begin{equation}{}\label{2-5}
\Psi_{{\mathcal L}} (f)\equiv\int_{l^1_+(X)}\exp(-\sum_k r_k f(x_k))
\,d{\mathcal L}({\xi})=\exp(-\int_X \log f(x)\,dm(x)),
\end{equation}
where $f$ is an arbitrary nonnegative measurable function on
$(X,m)$ with $\int\log f(x)dm(x)<\infty$. Note that, like the classical
Laplace transform, equation~\eqref{2-5} makes sense also for
complex-valued functions $f$ with nonnegative real component.

The construction of the integral model essentially uses only two properties of the
measure ${\mathcal L}$: its projective invariance with respect to
the group of multiplicators $M_a$ (namely, $d{\mathcal L}(M_a{\xi})=\phi(a)\,d{\mathcal L}({\xi} )$, where
$\phi(a)=  e^{\frac 12 \int_X\log a(x)d\,m(x)}$) and its ergodicity with respect
to this group. The invariance follows directly from the formula for the Laplace
transform; for a proof of the ergodicity, see \cite{V1}. In particular, the projective
invariance of ${\mathcal L}$ implies the orthogonality (respectively,
unitarity) of the representation
$\opn{INT}\,T$, and the ergodicity of ${\mathcal L}$
implies the irreducibility of $\opn{INT}\,T$.
For a detailed study of the properties of this remarkable measure, see
\cite{Vershik,V1}.

Further, in the construction of the integral model, the spaces
$H_r$, $r\in{\mathbb R {\,}}^*_+$ (which are the spaces of representations
of the subgroup $P_0$) are replaced by the countable tensor products
$H_{{\xi}}= \otimes _{k=1}^\infty H_{r_k}$,
${\xi}\in l^1_+(X)$, with stabilizing almost invariant vectors
$h\in H$ (which are the spaces of representations
$T_{{\xi}}$ of the group ${P}_0^X$).

By definition, the countable tensor product
$H_{{\xi}}=
\otimes _{k=1}^\infty H_{r_k}$ with stabilizing almost invariant vector
$h\in H$ is the completion of the inductive limit of the spaces
$\otimes_{k=1}^n H_{r_k}$ with respect to the embeddings
$x\in \otimes_{k=1}^n H_{r_k}
\mapsto x \otimes h_{r_{n+1}}\in \otimes_{k=1}^{n+1} H_{r_k}$, $h_{r_{n+1}}=h$.
The action of the operators corresponding to
the subgroup ${P}_0^X$ on $H_{{\xi}}$ is given
by the formula
$$
T_{{\xi}}(g(\cdot ))(\otimes _{k=1}^\infty f_k )= \otimes _{k=1}^\infty
(T_{r_k}(g(x_k))f_k).
$$
Note that the representations $ T_{{\xi}}$ are well defined only
for canonical representations $T$ of the subgroup ${P}_0$.

The space $\opn{INT}\,H$ of the integral model $U=\opn{INT}\,T$
is defined, by analogy with the space
$\mathcal H=\int_0^\infty H_r\,r^{-1}\,dr$ of the representation of
$P$ associated with a canonical representation $T$ of
${P}_0$, as the direct integral
of the Hilbert spaces $H_{{\xi}}$ with respect to the measure
${\mathcal L}$:
$$
\opn{INT}\,H = \int^{\oplus }_{l^1_+(X)}  H_{{\xi}}\,d{\mathcal L}({\xi}), \quad
H_{{\xi}} = \otimes_{k=1}^\infty H_{r_k},
$$
i.e., as the space of sections $F({\xi})=F(\{r_k,x_k\})$ of the fiber bundle over
$l^1_+(X)$ with fiber $ H_{{\xi}}$ over a point ${\xi}=\{r_k,x_k\}$.
The action of the subgroup  $P_0^X$  on the spaces $H_{\xi}$ induces a representation
$U=\opn{INT}\,T$ of this group in the whole space  $\opn{INT}\,H$:
$$
\opn{INT}\,T=\int_{l^1_+(X)}^ \oplus  T_{\xi}\,d{\mathcal L}(\xi),
\quad\text{where}\quad
T_{\xi}= \otimes _{k=1}^\infty T_{r_k,x_k};
$$
in more detail: $(U(g(\cdot ))F)({\xi})=\widetilde T_{{\xi}}(g(\cdot )\,
F({\xi})$. The extension of this representation of
$P_0^X$ to the whole group
$P^X$ is given by the formula
$$
(U(r_0(\cdot ))F)({\xi})=e^{\frac 12 \int_X\log r_0(x)\,dm(x)}\,
F(r_0(\cdot ){\xi})\quad\text{for}\quad r_0\in ({\mathbb R {\,}}^*_+)^X.
$$
It follows from the definition that the representation
$\opn{INT}\,T$ of the group ${P}^X$ is invariant under
transformations of the space $X$ preserving the measure $m$.

\begin{THM}{}\label{THM:21-80}(\cite{VG-7}) If a canonical representation
$T$ of the subgroup ${P}_0$ is irreducible, then the corresponding
representations
$\widetilde T$ and  $U{=}\opn{INT}\,T$ of the groups ${P}$ and
${P}^X$, respectively, are irreducible.
\end{THM}

\begin{REM*}{} The construction of the integral model also makes sense
if we replace the probability measure $m$ on $X$ with any positive finite measure
(with $m(X)= \theta $) or, which is equivalent, introduce the coefficient
$\theta>0 $ into the exponent in the definition of the Laplace transform
of the measure ${\mathcal L}$. The one-parameter family, thus defined, of representations
of ${P}^X$ still satisfies Theorem~\ref{THM:21-80}.
For simplicity, here we restrict ourselves with the case
$\theta =1$.
\end{REM*}

Let us define the spherical function of the representation
$U=\opn{INT}\, T$ of the group
${P}^X$ associated with a canonical representation $T$ of the subgroup
${P}_0$ in a space $H$ as the following function on ${P}^X$:
$$
\Psi(g)={\langle U(g)\Omega,\Omega \rangle},
$$
where $\Omega $ is a unit vector in the space $\widetilde H=\opn{INT} H$
of the form
\begin{equation}{}\label{7-1}
\Omega ({\xi})= \otimes _{k=1}^\infty (e^{-\frac 12 r_k}h_{r_k}),  \quad h_{r_k}=h
\quad\text {for}\quad {\xi}=\{r_k,x_k\}\in D_+(X).
\end{equation}

\begin{THM}{}\label{THM:1} The spherical function of the representation $U$
is given by the formula
\begin{equation}{}\label{4-4-4}
\Psi(g)=\exp\Bigl(\int_X(i\opn{Im} {\langle b(g(x)), f_0 \rangle}-
\frac 12 \|b(g(x)\|^2)\,dm(x)\Bigr),
\end{equation}
where $f_0(r)=e^{-\frac r2}\,h_r$ and
$b(g)= (\widetilde T(g) f_0)(r)-f_0(r)$ is the $1$-cocycle of the
representation
$\widetilde T$ of the group ${P}$ associated with $T$, see \eqref{14-4}.
\end{THM}

\section{Relation between integral models and Fock representations}

Given an arbitrary locally compact group $G$, each pair
$(\widetilde T, b)$, where $\widetilde T$ is a special representation of $G$
in a space $\mathcal H$ and $b:\,G\to\mathcal H$ is a nontrivial $1$-cocycle,
gives rise to a Fock model of the corresponding representation
of the current group
$G^X$, see \cite{VGG-6}. It is realized in the Hilbert space
$\opn{EXP}\mathcal H^X= \opl\oplus_{k=0}^\infty S^k \mathcal H^X$, where $S^k$
stands for the $k$th symmetric tensor power and
$\mathcal H^X=
\int_X^\oplus \mathcal H_x \,dm(x)$, $\mathcal H_x=\mathcal H$.
The operators of this representation are defined on the total subset of vectors
of the form $ \opn{EXP}v\,\,=\,\, \vac \,\,\oplus\,\, v \,\,\oplus\,\,
\frac{1}{\sqrt{2!}}\,v \otimes v \,\,\oplus
\,\,\frac{1}{\sqrt{3!}}\,v \otimes v\otimes v \,\, \oplus \,\,\cdots $
by the formula
$$
U(g)\opn{EXP}v=
\exp\Bigl(-\frac 12 \|b^X(g)\|^2-(\widetilde T^X(g)v,b^X(g))\Bigr)\,
\opn{EXP}(\widetilde T^X(g)v+b^X(g)).
$$
Here $\widetilde T^X$ and $b^X(g)$ are the representation of the
group $G^X$ in the space
$\mathcal H^X$ and the $1$-cocycle $b^X:\,G^X\to\mathcal H^X$ generated,
respectively, by the representation
$\widetilde T$ of the group $G$ and the $1$-cocycle $b:\,G\to\mathcal H$.
The operators $U(g_1g_2)$ and $U(g_1)U(g_2)$ in the Fock representation
are related by the formula $U(g_1g_2)=\exp(i\opn{Im}
\int_X a(g_1(x),g_2(x))\,dm(x))\,U(g_1)U(g_2)$,
where $a(g_1,g_2)=(\widetilde T(g_1)b(g_2),b(g_1))$. Thus, in the case of
unitary representations $\widetilde T$, the Fock representations are, in general,
projective. However, for groups of the form
${P}={\mathbb R {\,}}^*_+\rightthreetimes {P}_0$, the $2$-cocycle $a(g_1,g_2)$
is trivial, so that these representations are projective equivalent to
ordinary (nonprojective) representations.

\begin{THM}{}\label{THM:2}
Let $T$ be a canonical representation of the group
${P}_0$ and $\widetilde T$ be the special representation of the group $P$
associated with $T$ (see above). Then the integral model
$\opn{INT}\,T$  of the corresponding representation
of the group ${P}^X$ is projective
equivalent to the Fock representation of this group associated with
the same representation
$\widetilde T$ and the $1$-cocycle of the form \eqref{14-4} with
$f_0(r)=e^{-\frac r2}\,h_r$.
\end{THM}

\section{Integral models of representations of the subgroup $P \subset
SL(2,{\mathbb R {\,}})$ of triangular matrices and their extensions to the group
$SL(2,{\mathbb R {\,}})^X$}

Let us write elements of the subgroup $P  \subset SL(2,{\mathbb R {\,}})
$ of triangular matrices $g=\left(\begin{array}{ll}\epsilon^{-1}&0\\
\gamma & \epsilon\end{array}\right)$ as pairs
$(\epsilon, \gamma )$ with the multiplication law $(\epsilon _1, \gamma
_1)(\epsilon _2, \gamma _2)= (\epsilon _1\epsilon _2, \gamma
_1\epsilon_2^{-1}  +  \epsilon _1 \gamma _2)$. The group $P$
can be presented as the semidirect product
$P={\mathbb R {\,}}^*_+\rightthreetimes P_0$, where $P_0$
is the commutative group consisting of the elements $(\pm 1, \gamma )$.

Up to passing to conjugate representations,
the group $P_0$ has two
unitary irreducible canonical representations
$T^{\pm}$ in one-dimensional spaces $H^{\pm}$and one
orthogonal representation $T^0$ in a two-dimensional space $H^0$.
The representations $T_r^{\pm}$ in the spaces
$H^{\pm}_r=H^{\pm}$ are defined by the formulas
$T_r^{\pm}(1, \gamma )=e^{\pm ir^2 \gamma }\,\opn{id}$,
$T_r^{\pm}(\pm 1,0)=\opn{id}$. Accordingly, the countable tensor products
$ H^{\pm}_{{\xi}},$ ${\xi}=\{r_k,x_k\}$, in which the current group
$P_0^X$ acts, are also one-dimensional, so that the spaces
$\opn{INT}\,H^{\pm}$ of the representations
$U^{\pm}=\opn{INT}\,T^{\pm}$ of the group $P^X$
are the Hilbert spaces of complex-valued functionals
$F^{\pm}({\xi})=F(\{r_k,x_k\})$ on $l^1_+(X)$ with the norm
$$
\|F\|^2=\int_{l^1_+(X)} |F({\xi})|^2\,d{\mathcal L}({\xi}).
$$
The operators of these representations are identical on the center of $P^X$
and are uniquely
determined by the formulas
\begin{equation}{}\label{92-74}
(U^{\pm}(1, \gamma (\cdot ))F^{\pm})({\xi})=
e^{\pm i\sum r_k^2 \gamma (x_k)}\,F^{\pm}({\xi}),
\end{equation}
\begin{equation}{}\label{14-14}
(U^{\pm} (r_0(\cdot ),0)F)({\xi})=
e^{\frac 12 \int_X \log r_0(x)\,dm(x)}\,F^{\pm}( r_0 (\cdot )\,{\xi}).
\end{equation}

The orthogonal representations $T^0_r$ are realized in the two-dimensional
real subspaces $H_r^0 \subset H_r^+ \oplus H_r^-$ of vectors of the form
$(f,\overline f)$; the corresponding operators are the restrictions of the operators
$T^+_r(g) \oplus T^-_r(g)$ to these subspaces. Hence the spaces
$H^0_{{\xi}}$ of the representations of ${P}_0^X$
are countable tensor products of two-dimensional real spaces, and the integral
model $U^0=\opn{INT}T^0$ of the corresponding representation
of
${P}^X$ acts in the direct integral of these spaces with respect to the measure
${\mathcal L}$:
$\opn{INT}\,H^0=\int_{D_+(X)} H^0_{\xi}\,d{\mathcal L}(\xi)$.
The operators $U^0 (r_0(\cdot ),0)$ of this representation are given by the
same formula \eqref{14-14},
and the operators $U^0(1, \gamma (\cdot ))$ are given by the formula
\begin{equation}{}\label{192-74}
U^0(1, \gamma (\cdot )) (\otimes (f_k,\overline f_k))=
\otimes _{k=1}^\infty (e^{ir_k^2 \gamma (x_k)} f_k,
e^{-ir_k^2 \gamma (x_k)}\overline f_k).
\end{equation}
Note that under the natural embedding
$\opn{INT}\,H^0  \subset
\opn{INT}\,(H^+ \oplus H^-)$, the operators $U^0(g)$ are the restrictions
to $\opn{INT}\,H^0$ of the operators
$U(g)$ in the space $\opn{INT}\,(H^+ \oplus H^-)$.

Let us describe the extensions of the representations $U^{\pm}$ and $U^0$
of the current group $P^X$
to the group $SL(2,{\mathbb R {\,}})^X$.

Consider a homogeneous space of the group
$G=SL(2,{\mathbb R {\,}})$, namely, the upper complex half-plane $L$ endowed with
the following action of $G$:
$z\to  gz=\frac{\delta z+ \gamma }{\beta z + \alpha }$ for
$g= \begin{pmatrix}  \alpha & \beta \\ \gamma & \delta  \end {pmatrix}$
(the Lobachevsky plane). Denote by $L^X$ the space of bounded functions
$z: X\to L$. Note that the integral
$\int_X\log |z(x)|\,dm(x)$ converges for every function $z\in L^X$,
since $z$ is bounded.
The action of the group ${G}$ on $L$ induces the pointwise action of the group
${G}^X$ on $L^X$.

Further, consider on $L \times L$ the function
\begin{equation}{}\label{9-8}
c(z_1,z_2)=\log(-i(z_1-\overline{z_2} ))=\log[(v_1+v_2)-i(u_1-u_2)],
\quad z_k=u_k+iv_k,
\end{equation}
and define functions $u^{\pm}(g,z)$ on $G \times L$
by the formulas
\begin{equation}{}\label{9-7}
u^+(g,z)=\overline{u^-(g,z)}=
c(gz,gz_0)-c(z,z_0)-\frac 12 (c(gz_0,gz_0)-c(z_0,z_0)),
\end{equation}
where we have denoted $z_0=i$.

In order to extend the representations $U^{\pm}=\opn{INT}\,T^{\pm}$ of the group $P^X$
to the group $SL(2,{\mathbb R {\,}})^X$, with each function
$z\in L^X$ we associate the following functionals on $l^1_+(X)$:
$$
F^+_z({\xi})=e^{i\sum r_k^2z(x_k)}, \quad F^-_z({\xi})=\overline{F^+_z({\xi})}
=e^{-i\sum r_k^2\overline{z(x_k)}}.
$$
These functionals lie in the spaces
$\opn{INT}\, H^{\pm}$ and form total subsets in these spaces.

We define the action of the operators $U^{\pm}(g)$ of the
current group $SL(2,{\mathbb R {\,}})^X$
on the sets of functionals of the form $F^{\pm}_z$ by the following formula:
\begin{equation}{}\label{9-6}
U^{\pm}(g)\, F^{\pm}_z=e^{\frac 12 \int_X
u^{\pm}(g(x),z(x))\,dm(x)}\, F^{\pm}_{gz},
\end{equation}
where $u^{\pm}(g,z)$ is given by \eqref{9-7}.

\begin{THM}{}\label{THM:21} The operators $U^{\pm}(g)$
defined by \eqref {9-6} can be extended to
unitary operators on the corresponding spaces
$\opn{INT}\, H^{\pm}$ and determine
projective representations of the group
$G^X$ in $\opn{INT}\, H^{\pm}$. The restrictions of these
representations to the subgroup
$P^X$ coincide with the original representations of this subgroup.
\end{THM}

In order to extend the orthogonal representation
$U^0=\opn{INT}\,T^0$ of the group
${P}^X$ in the space $\opn{INT}\, H^0$
to the group $SL(2,{\mathbb R {\,}})^X$, with each function
$z\in L^X$ we associate the following functional on $l^1_+(X)$:
$$
F^0_z({\xi})= \otimes _{k=1}^\infty(2^{-1/2}\,e^{ir_k^2 z(x_k)},
2^{-1/2}\,e^{-ir_k^2 \overline{z(x_k)}}).
$$
The functionals $F^0_z$ lie in the space $\opn{INT}\, H^0$
and form a total set in this space.

We define the action of the operators $U^0(g)$ of the group ${P}^X$ on the set of
functionals of the form $F^0_z$ by the formula
\begin{equation}{}\label{9-10}
U^0(g)\, F^0_z=e^{\frac 12 \opn{Re}\int_X
u^{\pm}(g(x),z(x))\,dm(x)}\, F^0_{gz},
\end{equation}
where $u^{\pm}(g,z)$ is still given by \eqref{9-7}.

\begin{THM}{}\label{THM:22} The operators $U^0(g)$ can be extended to
orthogonal operators on $\opn{INT}\, H^0$ and determine
an orthogonal representation of the group
$G^X$ in $\opn{INT}\, H^0$. Its restriction to the subgroup
$P^X$ coincides with the original representation of $P^X$. The obtained representation
of the group
$SL(2,{\mathbb R {\,}})^X$ is equivalent to the representation constructed in
\cite{VGG-1}.
\end{THM}

\section{Extension of the integral models of representations of the maximal
parabolic subgroups of
$O(n,1), U(n,1)$ to the whole groups:
the orthogonal version}

It is known (see, e.g., \cite{VG}) that each of the groups $O(n,1)$ and $U(n,1)$
has a unique orthogonal irreducible special representation. Its
restriction to the maximal parabolic subgroup $P$ is also irreducible
and associated with a canonical orthogonal representation of the subgroup
$P_0$ appearing in the decomposition
$P={\mathbb R {\,}}^*_+\rightthreetimes P_0$; denote
this canonical orthogonal representation by $T^0$.
The question arises whether it is possible to extend the orthogonal representation
$\opn{INT}\,T^0$  of the group $P^X$
to a representation of the group $O(n,1)^X$ or $U(n,1)^X$, respectively.

\begin{THM}{}\label{THM:3} The integral models $\opn{INT}\,T^0$ of
orthogonal representations of the groups
$P^X$, where $P \subset O(n,1)$ or $P \subset U(n,1)$,
can be extended up to orthogonal representations of the groups
$O(n,1)^X$ and $U(n,1)^X$ that are equivalent to their Fock representations
described in \cite{VGG-6}. The corresponding intertwining operator
is generated by the mapping
$\Omega \mapsto \opn{EXP}0$ of the cyclic vectors.
\end{THM}

The assertion of the theorem follows from the coincidence of the spherical functions
of the corresponding representations of the group $P^X$.

In the case of the group $U(n,1)$, its special orthogonal irreducible representation
can be endowed with
the structure of a unitary representation in two ways.
The restrictions of the special unitary representations of
$U(n,1)$ thus obtained to the subgroup $P$ are irreducible and associated
with canonical unitary representations
$T^+$ and $T^-$ of the subgroup $P_0$.

\begin{THM}{}\label{THM:4} The integral models $\opn{INT}\,T^{\pm}$
of unitary representations of the group $P^X$ associated with the canonical
unitary representations
$T^{\pm}$ of the subgroup $P_0$ can be extended up to projective unitary
representations of the current group $U(n,1)^X$
that are projective equivalent to its Fock projective unitary representations
described in \cite{GG}.
\end{THM}

The extensions thus obtained are new models of representations of the current groups
$O(n,1)^X$ and $U(n,1)^X $, in which the representation space is realized
as the direct integral  of $P_0^X$-invariant subspaces
with respect to the measure
${\mathcal L}$,  and elements from
$({\mathbb R {\,}}^*_+)^ X \subset P^X$ permute these subspaces.

\section{Extension of the integral models of representations of the maximal
parabolic subgroup of $U(n,1)$ to the whole group: the unitary version}

Let us realize $U(n,1)$ as the group of linear transformations in
${\mathbb C {\,}}^{n+1}$ preserving the Hermitian form
$x_1\overline x_{n+1}+x_{n+1}\overline x_1+ |x_2|^2+ \ldots +|x_n|^2$,
and let us represent elements
$g\in U(n,1)$ in a block form:
$g=\|g_{ij}\|_{i,j=1,2,3}$, where the matrices on the diagonal are of orders
$1$, $n-1$, and $1$, respectively. In this realization, $P$ is the
subgroup in
$U(n,1)$ consisting of all lower triangular block matrices. It can be presented
as the semidirect product
$P={\mathbb R {\,}}^*_+\rightthreetimes P_0$, where
$P_0=D_0\rightthreetimes N$, $N$ is the Heisenberg group of order
$2n-1$ realized as the group of pairs $(t,z)\in{\mathbb R {\,}} \otimes
{\mathbb C {\,}}^{n-1}$, and $D_0\cong U(1) \times U(n-1)$.

In the construction of the integral models
$\opn{INT}\,T^{\pm}$ and
$\opn{INT}\,T^0$, the canonical unitary representations
$T^{\pm}$ of the group
$P_0$ are realized in the spaces $H^{\pm}$, respectively, of entire analytical
and entire anti-analytical functions
$f(z)$ on ${\mathbb C {\,}}^{n-1}$ with the norm
$\|f\|^2=\int_{{\mathbb C {\,}}^{n-1}} |f(z)|^2 e^{-zz^*}\,d\mu(z)$
(for the action of $P_0$ on $H^{\pm}$, see \cite{VG}).
The orthogonal canonical representation $T^0$ is realized in the subspace
$H^0 \subset H^+ \oplus H^-$ of vectors of the form $(f,\overline f)$.

According to the general definition, the representations
$U^{\pm}=\opn{UNT}\,T^{\pm}$ and $U^0=\opn{UNT}\,T^0$ of the group $P^X$
are realized in the space
$\opn{INT}\,H^{\pm}=\int^ \oplus _{D_+(X)}H^{\pm}_{{\xi}}\,d{\mathcal L}({\xi})$, where
$H^{\pm}_{{\xi}} = \otimes_{k=1}^\infty  H^{\pm}_{r_k}$, and in the space
$\opn{INT}\, H^0=\int^\oplus _{D_+(X)}H^0_{{\xi}}\,d{\mathcal L}({\xi})$, where
$H^0_{{\xi}} = \otimes_{k=1}^\infty  H^{\pm}_{r_k}$, respectively.

The extensions of these representations to the whole group
$U(n,1)^X$ can be constructed in a way similar to the case of
$SL(2,{\mathbb R {\,}})$. Namely, instead of the Lobachevsky plane,
we consider the homogeneous space $L$ of the group
$U(n,1)$ that is equivalent to the disk in ${\mathbb C {\,}}^n$ (the
$n$-dimensional complex Lobachevsky space):
$$
L=\{v=(a,b)\in{\mathbb C {\,}} \oplus {\mathbb C {\,}}^{n-1} \mid a+\overline a + b^* b< 0\},
\quad\text{ЗДЕ}\quad  b^* b=\sum \overline b_ib_i
$$
($b$ is a column vector). The action of the group $U(n,1)$ on $L$ is given by
the formula $g(a,b)=(a',b')$, where
\begin{gather*}{}
a'=(g_{11}+g_{12} b+g_{13} a)^{-1}\,(g_{31}+g_{32} b+g_{33} a),
\\
b'=(g_{11}+g_{12} b+g_{13} a)^{-1}\,(g_{21}+g_{22} b+g_{23} a).
\end{gather*}

Consider the following function on $L \times L$:
$$
c(v_1,v_2)=\log(-a_1-\overline a_2-b_2^* b_1)\quad\text{for}\quad  v_i=(a_i,b_i)\in L.
$$

Further, denote by $L^X$ the set of measurable bounded mappings
$v:\, X\to L$ with the pointwise action
$v \mapsto gv$ of the group
${P}^X$. Note that the integral
$\int_X c(v(x),v(x))\,dm(x)$ converges for every $v\in L^X$
 since  $v\in L^X$ is bounded.

With each triple $v=(a,b)\in L^X$, $r\in{\mathbb R {\,}}^*_+$, and $x\in X$,
we associate a vector $f_{v,r,x}\in H_r$:
$$
f_{v,r,x}(z)=\exp(r^2a(x)+r(z,b(x))), \quad  (z,b)=\sum z_i b_i;
$$
with elements $v\in L^X$ we associate the following functionals
$F_v^{\pm}({\xi})$ and $F_v^0({\xi})$ on $l^1_+(X)$:
$$
F^+_v({\xi})=\overline{F^-_v({\xi})}= \otimes _{k=1}^\infty f_{v,r_k,x_k}, \quad
F^0_v({\xi})= \otimes _{k=1}^\infty (2^{-1/2}f_{v,r_k,x_k},
2^{-1/2}\overline{f_{v,r_k,x_k}}).
$$
These functionals lie in the spaces
$\opn{INT}\, H^{\pm}$ and
$\opn{INT}\, H^0$, respectively, and form total subsets in these spaces.

On these sets we define the actions of the operators
$U^{\pm}(g)$ and $U^0(g)$
by the same formulas \eqref{9-6} and \eqref{9-10}, as in the case of the group
$SL(2,{\mathbb R {\,}})$, in which
$z$ should be replaced by $v$,
the function $c(z_1,z_2)$
should be replaced by the function $c(v_1,v_2)$ introduced above, and
instead of $z_0=i$  one should write
$v_0=(-1,0)$.

\begin{THM}{}\label{THM:90-1} The operators $U^{\pm}(g)$ and $U^0(g)$
can be extended to unitary and orthogonal operators in the spaces
$\opn{INT}\,H^{\pm}$ and $\opn{INT}\,H^0$, respectively;
they determine extensions of the original representations
of the  group $P^X$ up to projective unitary and orthogonal
representations of the group $U(n,1)^X$, respectively.
\end{THM}

\medskip
Translated by N.~V.~Tsilevich.


\begin{thebibliography}{99}

\bibitem{Araki} H.~Araki. Factorisable representations of the current algebra.
{\it Publ. RIMS Kyoto Univ.}, Ser.~A5, No.~3 (1970), 361--422.

\bibitem{VGG-1}
A.~M.~Vershik, I.~M.~Gelfand, and M.~I.~Graev.
Representations of the group $SL(2,R)$, where $R$ is a ring of
  functions.
  {\it Uspekhi Mat. Nauk} {\bf 28}, No.~5 (1973), 83--128.
  English translation in: {\it Representation Theory}, London Math.
  Soc. Lect. Note Ser. {\bf69}, Cambridge Univ. Press, 1982, pp.~15--60.

\bibitem{VGG-2}
A.~M.~Vershik, I.~M.~Gelfand, and M.~I.~Graev.
Irreducible representations of the group $G^X$ and cohomology.
{\it Funct. Anal. Appl.} {\bf 8} (1974), 151--153.

\bibitem{B}  F.~A.~Berezin. Representations of the continuous
direct product of universal coverings of the
group of motions of a complex ball.
{\it Trans. Mosc. Math. Soc.} {\bf 36} (1979), 281--298.

\bibitem{VGG-6}I.~M.~Gelfand, M.~I.~Graev, and A.~M.~Vershik.
Models of representations of current groups. In:
{\it Representations of Lie Groups and Lie
Algebras} (A.~A.~Kirillov, ed.), Akad\'emiai Kiad\'o, Budapest, 1985,
    pp.~121--179.

\bibitem{GG}I.~M.~Gelfand and M.~I.~Graev. Special representations of the
group ${\rm SU}(n,1)$ and projective unitary representations of the current
group ${\rm SU}(n,1)^X$. {\it Russian Acad. Sci. Dokl. Math.} {\bf 48},
No.~2  (1994), 291--295.

\bibitem{VG1} A.~M.~Vershik and M.~I.~Graev.
A commutative model of a representation of the group ${\rm O}(n,1)^X$
and a generalized Lebesgue measure in a distribution space.
{\it Funct. Anal. Appl.}  {\bf39},  No.~2  (2005), 81--90.

\bibitem{VG}  A.~M.~Vershik and M.~I.~Graev.
Structure of the complementary series and special representations
of the groups $O(n,1)$ and $U(n,1)$. {\it Russian Math. Surveys}
{\bf61}, No.~5 (2006), 799--884.

\bibitem{VG-indag}  M.~I.~Graev and A.~M.~Vershik.
The basic representation of the current group $O(n,1)^X $ in the
$L^2$ space over the generalized Lebesgue measure. {\it Indag. Math.} {\bf16},
No.~3/4  (2005), 499--529.

\bibitem{VG-7}   A.~M.~Vershik and M.~I.~Graev.
Integral models of representations of current groups.
{\it Funkts. Anal. i Prilozh.} {\bf 42}, No.~1 (2008), 22--32.

\bibitem{Vershik}
 A.~M.~Vershik. Does there exist the Lebesgue
 measure in the infinite-dimensional space?
 {\it Proc. Steklov Inst. Math.} {\bf259} (2007), 256--281.

\bibitem{V1} A.~Vershik. Invariant measures for continual Cartan
subgroup. {\it J. Funct. Anal.} (2008).

\bibitem{VK}  A.~M.~Vershik and S.~I.~Karpushev.
Cohomology of groups in unitary representations, the neighborhood of the identity,
 and conditionally positive definite functions.
{\it Math. in USSR}  {\bf47} (1984), 513--526.

\bibitem{VT} A.~M.~Vershik and N.~V.~Tsilevich. Fock factorizations
 and decompositions of the $L^2$ spaces over general Levy processes.
 {\it Russian Math. Surveys}  {\bf58}, No.~3 (2003), 427--472.

\bibitem{TVY} N.~Tsilevich, A.~Vershik, and M.~Yor.
An infinite-dimensional analogue of the Lebesgue measure and distinguished
properties  of the gamma process. {\it J. Funct. Anal.} {\bf185}, No.~1
(2001), 274--296.

%\bibitem{W} J.~A.~Wolf. {\it Unitary Representations of Maximal Parabolic Subgroups
%of the Classical Groups}. Mem. Amer. Math. Soc. {\bf8}, No.~180 (1976).

\bibitem{Is} R.~S.~Ismagilov. {\it Representations of Infinite-Dimensional Groups}.
Transl. Math. Monographs {\bf152}. Amer. Math. Soc., 1996.

%\bibitem{K}A.~A.~Kirillov. Unitary representations of nilpotent Lie groups.
%{\it Uspekhi Mat. Nauk}  {\bf17}, No.~4(106) (1962), 57--110.

%\bibitem{K-2002}
%A.~A.~Kirillov. {\it Lectures on the Orbit Method}. Novosibirsk, 2002.

%\bibitem{Perelomov} A.~Perelomov. {\it Generalized Coherent States and
%Their Applications}.
%Springer-Verlag, Berlin, 1986.

%\bibitem{Hurt}N.~E.~Hurt. {\it Geometric Quantization in Action}.
%D.~Reidel Publishing Co., Dordrecht--Boston, Mass., 1983.

%\bibitem{BE} H.~Bateman and A.~Erdelyi.
%{\it Higher Transcendental Functions}, Vol.~2.
Moscow, 1974.

\bibitem{TsV} B.~Tsirelson and  A.~Vershik. Examples of nonlinear continuous tensor product of measure
spaces and non-Fock factorizations. {\it Rev. Math. Phys.} {\bf10}, No.~1 (1998),
81--145.

\bibitem{XR}
E.~Hewitt and K.~A.~Ross. {\it Abstract Harmonic Analysis},
Vol.~II.® Springer-Verlag, New York--Berlin, 1970.

\end{thebibliography}
\end{document}